\RequirePackage{fix-cm}

\documentclass{article}

\usepackage{amsmath,amsfonts,amsthm}
\usepackage[utf8]{inputenc} %unicode support

\usepackage{graphicx}
\newcommand{\Figure}{\textsc{Fig.}}
\newcommand{\bdelta}{\boldsymbol{\delta}}
\theoremstyle{remark}
\newtheorem*{remark}{Remark}
\newcommand{\V}[1]{\underline{#1}}                     % vectores

\newcommand{\iC}{\ensuremath{\widetilde{C}}}
\newcommand{\iI}{\ensuremath{\widetilde{I}}}

\title{On the computation of plate assemblies using realistic 3D joint model: a non-intrusive approach}

\author{
{Guillaume Guguin}$^1$, {Olivier Allix},$^1$ {Pierre Gosselet}$^1$, {Stéphane Guinard}$^2$ \\
  $1$ {LMT-Cachan, ENS-Cachan/CNRS/Universit\'e Paris-Saclay},\\ 
  {61 avenue du Pr\'esident Wilson}, {94235} {Cachan}, {France} \\                                   % country
  $2$ {Airbus Group Innovation},\\
  {18 rue Marius Tercé Zac St Martin du Touch}, {31000} {Toulouse},  {France}
}

\begin{document}
\maketitle
%%%%%%%%%%%%%%%%%%%%%%%%%%%%%%%%%%%%%%%%%%%%%%
%%                                          %%
%% The Abstract begins here                 %%
%%                                          %%
%% Please refer to the Instructions for     %%
%% authors on http://www.biomedcentral.com  %%
%% and include the section headings         %%
%% accordingly for your article type.       %%
%%                                          %%
%%%%%%%%%%%%%%%%%%%%%%%%%%%%%%%%%%%%%%%%%%%%%%

\begin{abstract} % abstract
%\parttitle{First part title} %if any
%Text for this section.
%
%\parttitle{Second part title} %if any
%Text for this section.

Most large engineering structures are described as assemblies of plates and shells and they are computed as such using \textit{adhoc} Finite Element packages. In  fact their computation in 3D would be much too costly.  In this framework, the connections between the parts are often modeled by means of simplified tying models. In order to improve the reliability of such simulations, we propose to apply a non-intrusive technique so as to virtually substitute the simplified connectors by a precise 3D nonlinear model, without modifying the global plate model. Moreover each computation can be conducted on independent optimized software.
%To overcome this limitation, we propose to extend  a non-intrusive coupling techniques between 2D and 3D models to the nonlinear case using dedicated software for the 2D and 3D computations. 
%The principle is, through iterations, to substitute locally the simplified connector by a 3D nonlinear model of the junction while keeping the global 2D model unchanged. 
After a description of the method, examples are used to analyze its performance, and to draw some conclusions on the validity and limitation of both the modeling of junction by rigid connectors and the use of submodeling techniques for the estimation of the carrying capacity of bolted plates. 

\textbf{keywords: assembly, non-intrusive coupling, bolt}
\end{abstract}

%%%%%%%%%%%%%%%%%%%%%%%%%%%%%%%%%%%%%%%%%%%%%%
%%                                          %%
%% The keywords begin here                  %%
%%                                          %%
%% Put each keyword in separate \kwd{}.     %%
%%                                          %%
%%%%%%%%%%%%%%%%%%%%%%%%%%%%%%%%%%%%%%%%%%%%%%

%%%%%%%%%%%%%%%%%%%%%%%%%%%%%%%%%%%%%%%%%%%%%%
%%                                          %%
%% The Main Body begins here                %%
%%                                          %%
%% Please refer to the instructions for     %%
%% authors on:                              %%
%% http://www.biomedcentral.com/info/authors%%
%% and include the section headings         %%
%% accordingly for your article type.       %%
%%                                          %%
%% See the Results and Discussion section   %%
%% for details on how to create sub-sections%%
%%                                          %%
%% use \cite{...} to cite references        %%
%%  \cite{koon} and                         %%
%%  \cite{oreg,khar,zvai,xjon,schn,pond}    %%
%%  \nocite{smith,marg,hunn,advi,koha,mouse}%%
%%                                          %%
%%%%%%%%%%%%%%%%%%%%%%%%%%%%%%%%%%%%%%%%%%%%%%

%%%%%%%%%%%%%%%%%%%%%%%%% start of article main body
% <put your article body there>

\section{Introduction}

The simulation of large structures undergoing complex local nonlinear  phenomena is still a major scientific and industrial challenge. One of the main difficulties originates from the difference of length scale between the global response of the structure and the localized phenomena. 
To address those problems a first type of computational approach is based on homogenization, as FE\textsuperscript{2} \cite{feyel2000fe} but it works well as long as the scales are sufficiently separated. 
To overcome this limitation concurrent multiscale methods have been developed. They are often based on domain decomposition techniques like FETI \cite{wyart_substructured_2007}, FETI-DP\cite{amini2009multi} or the LATIN multiscale method \cite{Ladeveze00a, Ladeveze:2001ab} and its optimization for the multiscale treatment of nonlinear problems \cite{Kerfriden09,Saavedra11}.\medskip

Moreover most of large industrial structures  are described as an assembly of plates and shells, whereas local phenomena often require 3D models to be properly analyzed. To deal with such problems, several methods have been applied or developed for the coupling of 2D and 3D models, like the Arlequin method \cite{dhia1998problemes,dhia2005arlequin}, transition elements \cite{gmur1993three,garusi_hybrid_2002}, MPCs approaches \cite{mccune_mixed-dimensional_2000,donaghy_dimensional_1996} or Nitsche's method \cite{nguyen2014nitsche}.

Most of these methods are quite demanding in terms of software development and therefore they are seldom used in industrial packages. To overcome these drawbacks, non-intrusive approaches have recently been proposed \cite{gendre2009non}. They are nowadays the subject of extensions and developments: thermoelasticity with GFEM/FEM coupling \cite{plews_improved_2011}, crack propagation in XFEM/FEM coupling \cite{duval2014non},  stochastic simulations \cite{safatly2012methode} and dynamics \cite{bettinotti2013coupling,bettinotti_allix_perego_oancea_malherbe_2014}. 

In \cite{guguin2014nonintrusive} a non-intrusive coupling between plate and 3D models was proposed in the case of linear behaviors. 
The present paper concerns the extension of this approach to the simulation of bolted assemblies of plates where bolts are described with full 3D nonlinear models.  Such structures are good candidates for the iterative global-local non-intrusive strategy for two main reasons. First, the  detailed computation of a tightened bolt with frictional contact on all surfaces is a very hard and time-consuming task to perform using commercial software. Second, the construction of the reference model, which would correspond to the assembly of a plate model and a 3D model for the bolt, would be very complex. The non-intrusive framework provides answers to both these problems. First, it allows the use of dedicated software for the local computations (in our case, COFAST a parallel software based on the LATIN domain decomposition method \cite{Champaney1999249}). Second, it allows the easy coupling of a general Finite Element software (here Code\_Aster from EDF) for the plate computation, with COFAST, because the global model is unchanged during the iterative process. These properties were exploited in \cite{daghia2012micro} for the simulation of damage in composite laminates at the meso and the micro scales using dedicated pieces of software.\medskip

The non-intrusive framework aims at solving the reference problem iteratively, by solving at each iteration both the global problem with prescribed residual traction at the interface and the local nonlinear problems submitted to prescribed displacement. Several techniques have been proposed to improve the convergence rate of the method by means of acceleration techniques \cite{gendre2009non, duval2014non, liu2014} or improved interface conditions \cite{gendre2011two}.
The method has therefore common points with so called nonlinear domain decomposition methods (or nonlinear relocalization techniques) \cite{klawonn14,NEGRELLO.2016.1} which proved their efficiency and gain in robustness in the case of buckling \cite{Cresta07}, post-buckling \cite{Hinojosa2014} and damage analysis \cite{Bordeu2009,allix2010}. Other proposals have been made to take into account the fact that the phenomena of interest are localized,  aiming at a better representation of the target model  \cite{duster2007,Duster20073524,plews_improved_2011}.\medskip

The paper is organized as follows. In section~\ref{sec:presNI}, a summary of the reference problem corresponding to the coupling of 2D and 3D models according to \cite{guguin2014nonintrusive} is presented. The non-intrusive algorithm is presented in section~\ref{sec:NIalgorithm}. 
In section~\ref{sec:glob-local-res}, the results of the iterative coupling are analyzed in the case of a bolted joint.  These results are compared to those corresponding to the plate solution and to the submodeling technique. In section~\ref{sec:params}, the influence of some parameters is assessed regarding the rate of convergence of the iterative process and regarding the accuracy of the coupled model compared to a full 3D solution.

%%%%%%%%%%%%%%%%%%%%%%%%%%%%%%%%%%%%%%%%%%%%%%%%%%%%%%%%%%%%%%%%%
\section{The reference problem}\label{sec:presNI}

\begin{figure}[ht]
\centering
\includegraphics[width=0.9\linewidth]{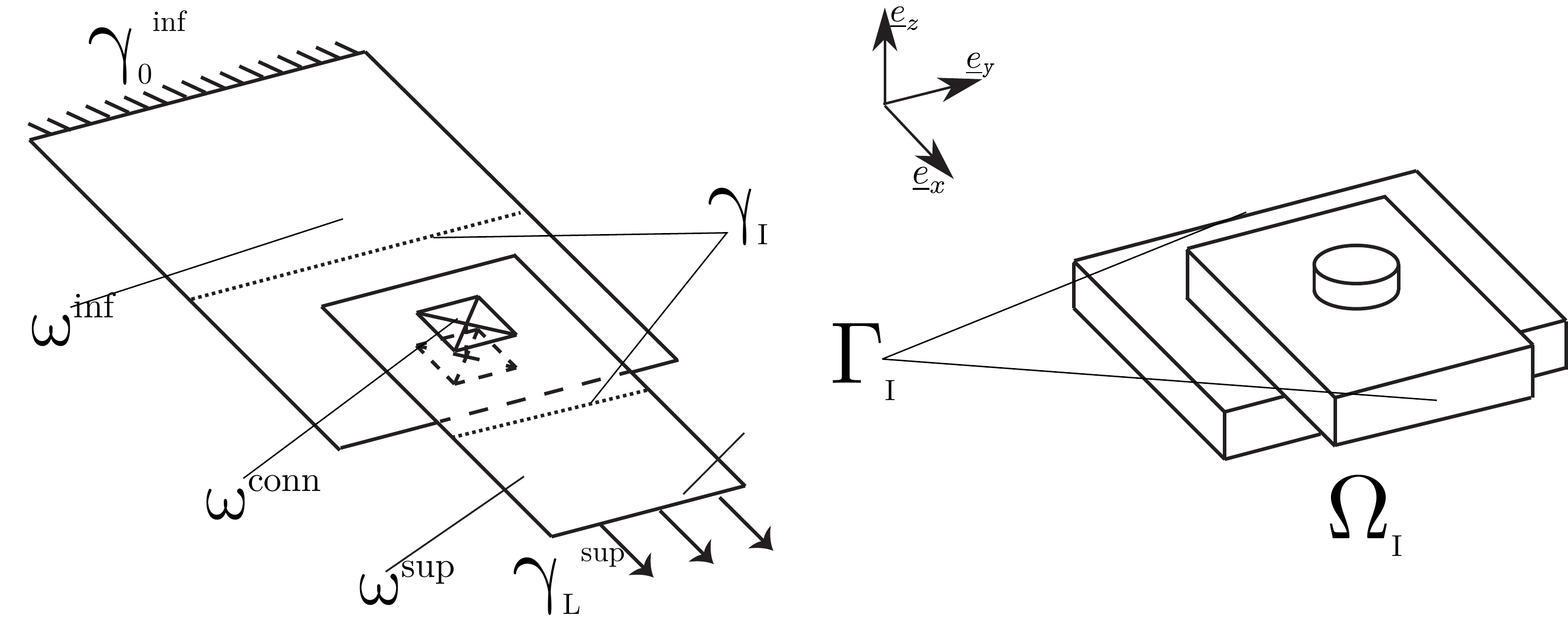}
\caption{Left, global model: two plates assembled by one 1D connector. Right, local model: full representation of the bolt with surrounding plates.}\label{fig:01}
\end{figure}

We consider the typical problem of two plates connected by a bolt. Starting from a simplified model of the assembly using plates elements and a simple connector (typically a beam), our aim is to perform the non-intrusive substitution of the connector by a full 3D model of the bolt. A sketch of the two models is presented in \Figure~\ref{fig:01}. 
In what follows lower case Greek characters are used to describe the geometry of the plates whereas  capital Greek characters are used for the geometry of the 3D domains:
any 3D domain $\Omega_X$ possesses a plate counterpart  $\omega_X$, so does any interface $\Gamma_Y$ whose plate counterpart is written $\gamma_Y$.

\begin{figure}
\includegraphics[width=0.9\linewidth]{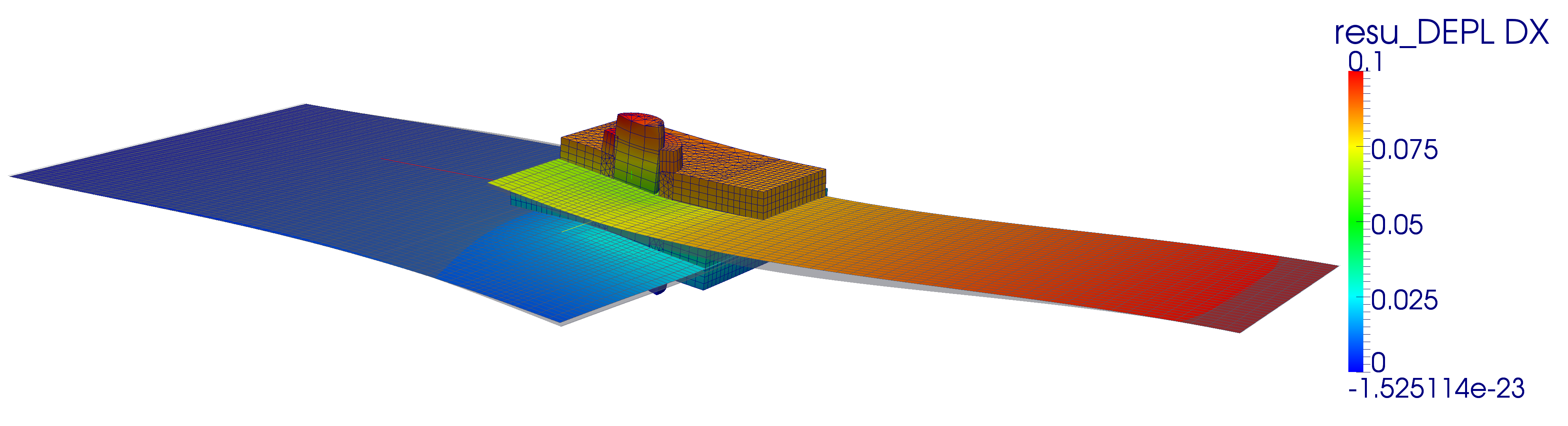}
\caption{Deformed shape of the targeted coupled problem.}\label{fig:02}
\end{figure}

The target hybrid model, which will be precisely defined in this section, corresponds to plate models connected by a full 3D model of the bolt. At the convergence of the iterations, the solution of the hybrid model is obtained, as illustrated in \Figure~\ref{fig:02}. Let us note that since the problem is symmetric (in the $(\underline{e}_x,\underline{e}_z)$ plane) only half of the bolt is computed and shown on the figure. 

%Typically $\Omega = \omega\times\left[\frac{-h}{2},\frac{h}{2}\right]$, where $h$ is the thickness of the plate.

\subsection{The plate model}

The global plate model is the assembly of two plates, $\omega^{inf}$ and $\omega^{sup}$ respectively the lower plate and the top plate, and a rigid connector between them $\omega^{conn}$ \Figure~\ref{fig:01} (left). The plates are $20$ mm thick and  $280$ mm long. The lower plate is $160$ mm wide whereas and the top plate is  $80$ mm wide. In order to simplify the presentation, the plates are assumed to be made out of homogeneous isotropic linear elastic material with Young modulus $E =200$ GPa and Poisson ratio $\nu = 0.3$; the handling of orthotropic composite plates is explained in \cite{guguin2014nonintrusive}.

The lower plate $\omega^{inf}$ is clamped on the part $\gamma^{inf}_{0}$ of its boundary and a prescribed tension displacement ($\underline{e}_x$ direction) is applied on the part $\gamma^{sup}_{L}$ of the top plate's boundary, \Figure~\ref{fig:01}:
\begin{equation}
\left\lbrace
\begin{aligned}
\V{U} &= \V{0} & \text{on}\; \gamma^{inf}_{0}\\
\V{U} &= U\,\V{e}_{x} & \text{on}\; \gamma^{sup}_{L}
\end{aligned}\right.
\end{equation}
Note that this tension is applied after preload (tightening) was applied to the bolt (see the description of the patch).

%To set up the plate problem, we will use the following notation: we call $\omega^{1}$ and $\omega^{2}$ respectively the lower and the top plate.

A classical Reissner-Mindlin plate formulation is used, and the contact is not taken into account at this stage. The kinematic is given by:
\begin{equation}\label{eq:platedep}
\V{U}(m,z)=\V{V}(m) + w(m)\V{e}_z - \V{e}_z\wedge\V{\theta} = \V{V}(m) + w(m)\V{z} - z
\begin{bmatrix}
\theta_{y}(m)\\ 
-\theta_{x}(m) \\
0
\end{bmatrix}
\end{equation}
where $m \in \omega^{inf} \cup \omega^{sup}$ is a point of the mid-surface of the plates and $\V{e}_z$ the normal direction (in the thickness). The stress in the plates is described with the eight generalized forces ($\underline{\underline{N}}$ tension, $\underline{\underline{M}}$ bending, $\underline{Q}$ shear): 
\begin{equation}\label{eq:defF}
 \mathbb{F} = \left( N_{xx},N_{yy},N_{xy},M_{xx},M_{yy},M_{xy},Q_{x},Q_{y} \right)
\end{equation} 

\begin{figure}[ht]
\centering
	\includegraphics[width=0.5\linewidth]{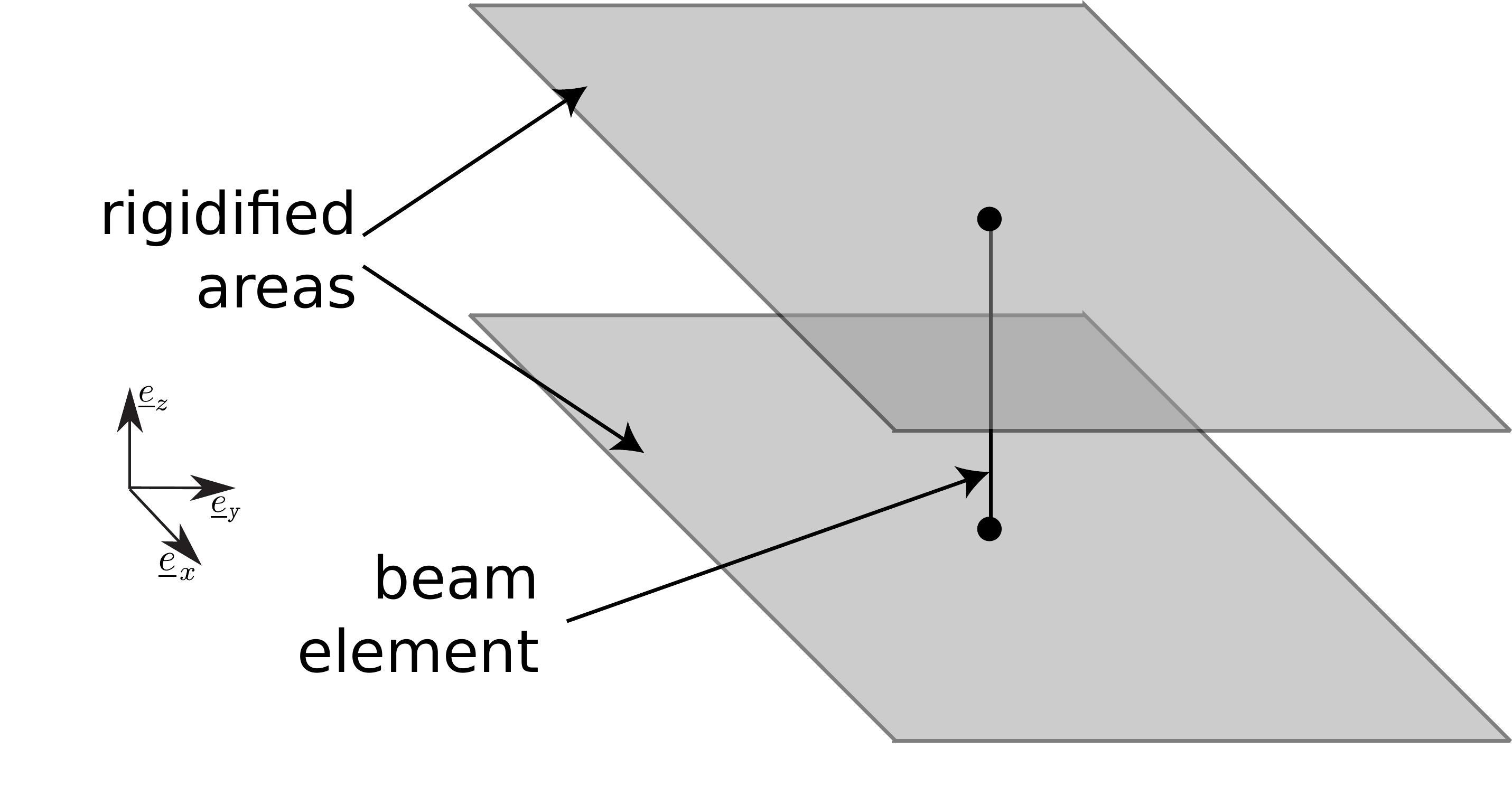}
\caption{The definition of the connector.}\label{fig:03}
\end{figure}

Since the connector will ultimately be virtually replaced by the 3D model described in next subsection, a very crude connector is chosen: a perfectly rigid beam element surrounded by a rigid region of the size of the bolt  ($20$ mm) in order to avoid localization around one node, see \Figure~\ref{fig:03}.
%
%This connector is sufficient to get an approximation global response. Actually, in the next section, this part will be \textit{virtually} replaced by the 3D model, so an accurate definition of the connector is not required. 
%
Beside its simplicity, this connector has the advantage to be very rigid compared to its 3D counterpart, which is a sufficient condition to ensure the convergence of the coupling algorithm \cite{gendre2009non,guguin14phd}.

%In addition, some part of this model will be used to compute the local model. So, the two sets of interfaces , $\gamma_I$ and $\gamma_C$ are defined in conjunction with the 3D model. The $\gamma^I$ interfaces define the strict domain where the solution is substituted and where the plate solution has no mecanical meaning (inside $\omega^I$). The $\gamma_C$ interfaces are define with the buffer zone, which is a small overlapping of the 3D model design to absorb the small perturbation after the application of the boundary conditions.

%an interface $\gamma_I$ is defined, as the 1D intersection between the 3D and the 2D model. This interface will be used to extract the global displacement of the solution.
%Another interface, $\gamma_C$ is set 

\subsection{The 3D patch}

\begin{figure}[ht]
\centering
\includegraphics[width=0.90\linewidth]{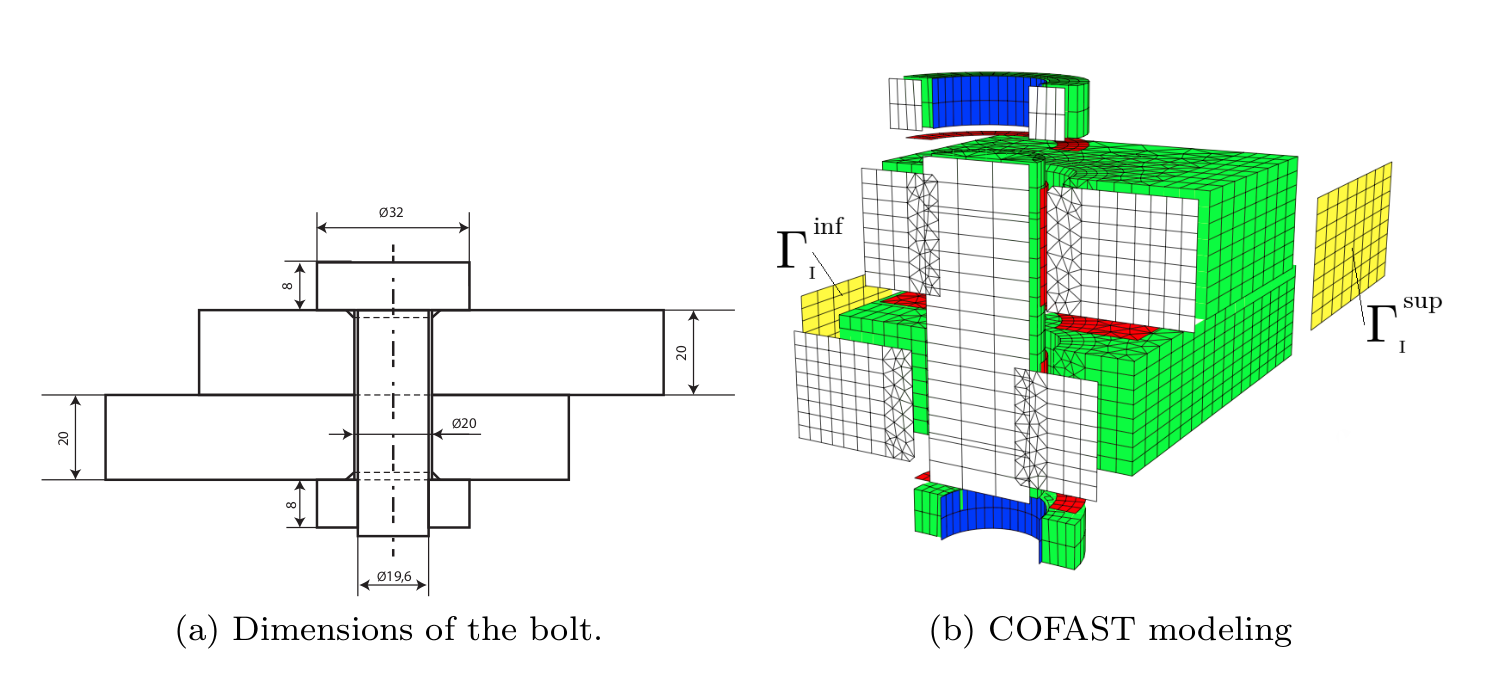}
%\subfloat[][Dimensions of the bolt.]{\label{fig:cotation_3D}
%\includegraphics[width=0.45\linewidth]{cotation_3D.pdf}
%}
%\subfloat[][COFAST modeling]{\label{fig:eclate}
%\includegraphics[width=0.45\linewidth]{eclate.png}
%}
\caption{3D representation of the bolt.}\label{fig:04}
\end{figure}

The local patch is designed to replace the connector. Thus, a full 3D representation of the bolt is used with unilateral frictional contact interfaces between each part of the assembly. The dimensions of the bolt are given in \Figure~\ref{fig:04}(a). The problem being symmetric, only half of the bolt is computed. The material used for the screw and the nut is linear elastic with Young's modulus $E=300$ GPa and Poisson's ratio $\nu=0.3$. Coulomb's friction coefficient is equal to $0.3$ on all interfaces. 

To solve this model, we use the dedicated contact solver COFAST, developed in the LMT-Cachan by L. Champaney \cite{Champaney1999249}. 
This software uses the LATIN method \cite{ladeveze1999nonlinear} where interfaces between subdomains are the support of the contact conditions, see \Figure~\ref{fig:04}(b).

%The nut have an internal diameter of $\varnothing 10mm$ and the diameter of the head is $\varnothing 32mm$. 

Note that before coupling with the global problem, preload is applied by enforcing relative displacements between the nut and the screw.  It is adjusted to match realistic values of tension in the screw ($200$ MPa).

\subsection{Connections between the models}

The plates+connector model describes the entirety of the structure and occupies the (2D+1D) domain $\omega$.
The 3D model of the bolt occupies what we call the zone of interest $\Omega_I$. In the original method \cite{gendre2009non}, the coarse modeling of $\omega_I=\omega\cap\Omega_I$ was simply replaced by $\Omega_I$ through iterations; in \cite{guguin2014nonintrusive}, it was shown that because of the edge effects which affect plate solutions, it is interesting to introduce a zone of transition between the two models. 

\begin{figure}[ht]
\centering
\includegraphics[width=0.9\linewidth]{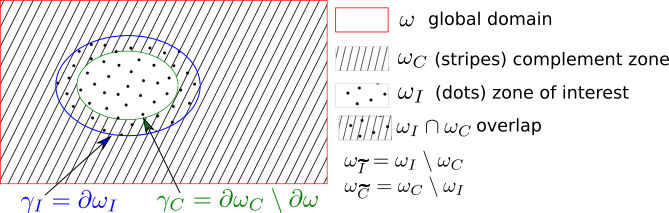}
\caption{The domain decomposition with overlap and the interfaces.}\label{fig:05}
\end{figure}

This idea is sketched in \Figure~\ref{fig:05}:
the 3D model is the only one taken into account in the inner zone of interest $\Omega_{\iI}\subset\subset\Omega_I$, the plate model is the only one taken into account in the outer complement zone $\omega_{\iC}\subset\subset\omega_C$, there exists an overlap $\Omega_I\cap\omega_C$ also called buffer zone where the two models are equivalent in a certain sense. 

\begin{figure}[ht]
\centering
\includegraphics[width=0.8\linewidth]{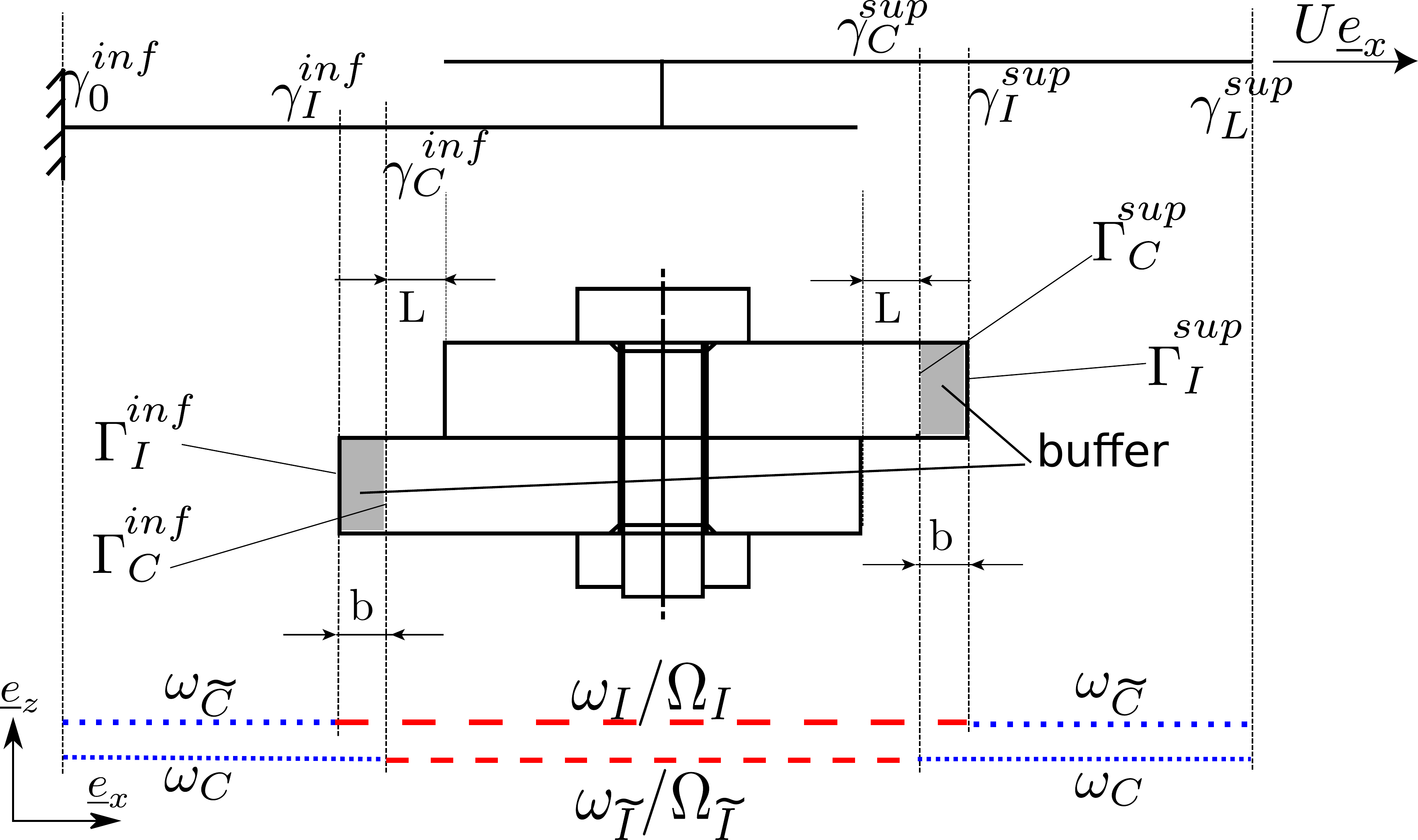}
\caption{Definition of the size of the local model and of the buffer zone.}\label{fig:06}
\end{figure}

%In \cite{guguin2014nonintrusive}, better results were obtained when a small buffer area was considered: the domains $\Omega^{I}$ and $\omega^{C}$ are defined with a small overlap, called the buffer zone. 
Two interfaces are defined: $\gamma_{I}=\partial\Omega_{I}\cap\omega_{C}$ and $\gamma_{C}=\partial\omega_{C}\cap\Omega_{I}$. 
In order to simplify the example, regular and ruled quadrilateral meshes are used along the interfaces. Thus, $\gamma_I$ and $\gamma_C$ are straight lines between two conforming meshes, see \Figure~\ref{fig:06}. 
 
The substituted solution is expressed with the following form:
\begin{equation}\label{eq:deplhybrid}
 \V{u}^{hyb} = \left\lbrace
 \begin{aligned}
& \V{u}^{L}   \text{ in } \Omega_{I} \text{ solves the 3D equations} \\
& \V{u}^{G}   \text{ in } \omega_{C} \text{ solves the plate equations}  \end{aligned}\right.
\end{equation} 
and it satisfies the following transmission conditions (in a sense which is made more precise afterwards):
\begin{equation}\label{eq:interfaces}
\left\lbrace \begin{aligned}
& \V{u}^{L} \text{ and } \V{u}^{G} \text{ are continous on }\Gamma_I \\
& \text{3D stress  and plate generalized forces are balanced on }\gamma_C \\
 \end{aligned}\right.
\end{equation} 

To be more specific, let us assume that a finite element discretization is used: $\mathbf{u}^G$ is the vector of plate unknowns defined on $\omega$, $\mathbf{u}^L$ is the vector of 3D unknowns defined on $\Omega_I$, $\mathbf{f}_{ext}$ stands for the generalized forces (they include the effect of imposed displacements), $\mathbf{K}$ for the stiffness matrices, and $\mathbf{f}_{int}$ for the internal forces (work of the stress field in the finite element subspace).\medskip

The global plate model is enriched by an extra loading $\bdelta$ which is only non-zero on the degrees of freedom of the interface $\gamma_C$ (the values of $\bdelta$ are determined by the coupling algorithm). The global plate model equilibrium thus can be written as:
\begin{equation}\label{eq:Gloeq}
\mathbf{K}^G \mathbf{u}^G = \mathbf{f}_{ext}^G + \boldsymbol{\delta}
\end{equation}
Note that the assumption of a linear global model simplifies the explanation and speeds up the solving (because the factorization is done only once) but it is not a requirement of the coupling algorithm.

The local 3D equilibrium submitted to boundary conditions inherited from the plate computation can be written as:
\begin{equation}\label{eq:Loceq}
\left\{\begin{aligned}&\mathbf{f}_{int}^L(\mathbf{u}^L)  = -\mathbf{f}_{ext}^L   \text{ in } \Omega_{I} \\
& \mathbf{u}^{L} = \mathbf{R}\mathbf{u}^{G} +\mathbf{w}^G(\mathbf{u}^{G}) \text{ on } \Gamma_{I} 
\end{aligned}\right.
\end{equation}
Matrix $\mathbf{R}$ represents the classical 3D rigid section movements. It enables to express the 3D Reissner-Mindlin kinematic from the plate degrees of freedom ($\mathbf{u}^G$).  Because the rigid section kinematic is a very coarse hypothesis, we supplement it with warping displacements $\mathbf{w}^G(\mathbf{u}^{G})$, proportional to the plate stress state. To be more precise, the applied displacement on $\Gamma_I$ has the following form:
\begin{equation*}
\V{u}^{L}(m,z)=\underset{\mathbf{R}\mathbf{u}^{G}}{\underbrace{\V{u}(m) + w(m)\V{e}_z - \V{e}_z\wedge\V{\theta}}}+ \underset{\mathbf{w}^G(\mathbf{u}^{G}) }{\underbrace{\sum_{i=1}^8 {\V{\mu}^{I}_i(z)\mathbb{F}_i^{\iC}}}}
\end{equation*}
with $(\mathbb{F}^{\iC}_i)$ the eight generalized forces extracted along the interface $\gamma_I$ \eqref{eq:defF}. $(\V{\mu}^I_i(z))$ is the basis of warping adapted to the stacking and materials which was determined before the computation by solving preliminary problems \cite{guguin2014nonintrusive}. 
These preliminary problems conducted on a 3D representative volume element of the composite plate enable us to compute distributions of stress and displacement (warping) associated with any solution of Saint-Venant problems. In turn, Saint-Venant problems are associated with specific components of the plate generalized stresses. 

From the global and local computations, it is possible to post-process the plate and 3D nodal reactions on the $\gamma_C/\Gamma_C$ interface:
\begin{equation}\label{eq:pplam}
\begin{aligned}
& \boldsymbol{\lambda}^{C} = \left(\mathbf{K}^{C}_{|\omega_C}\mathbf{u}^G - \mathbf{f}^{C}_{ext|\omega_C}\right)_{|\gamma_C} \\
&\boldsymbol{\lambda}^{L} = -\left(\mathbf{f}_{int|\Omega_{\iI}}^{L}(\mathbf{u}^{L})  +\mathbf{f}_{ext|\Omega_{\iI}}^{L}\right)_{|\Gamma_C}
\end{aligned}
\end{equation}

Finally the balance of the nodal reactions is enforced on the static plate quantities, thanks to the transpose operator $\mathbf{R}^T$ which enables to compute the generalized forces associated to the 3D nodal reactions $\boldsymbol{\lambda}^L$:
\begin{equation}\label{eq:balance}
\boldsymbol{\lambda}^{C} + \mathbf{R}^{T}\boldsymbol{\lambda}^{L} = 0  
\end{equation}

In spite of the care taken to recover the 3D boundary conditions on $\Gamma^I$ from plate quantities, small edge effects may occur, in particular in the case of stratified plates where suppressing edges effects require high order theories which are not always available in finite element analysis software or which could difficultly be implement in a non-intrusive manner. The  role of the buffer zone (the overlap between $\Omega^I$ and $\omega^C$) is then to dampen these edge effects and make the evaluation of the 3D stress reliable on $\Gamma^C$. In the example below, the size of the buffer zone	 is set to $b=5$mm, which was sufficient to absorb the local artificial effects. It corresponds to four elements between the two interfaces $\Gamma_{C}$ and $\Gamma_{I}$. The dimension $b$ is represented on  \Figure~\ref{fig:06}, together with dimension $L$ which characterizes the size of the 3D domain of interest $\Omega_I$. $L$ corresponds to the size of the part of the 3D domain which is not strictly necessary to represent the bolt correctly but which was inserted as a way to keep the 3D/2D transition away from the zone dominated to by 3D effects.  

\section{The non-intrusive iterative algorithm}~\label{sec:NIalgorithm}

The system (\ref{eq:Gloeq},\ref{eq:Loceq},\ref{eq:pplam},\ref{eq:balance}) can be interpreted as finding the traction $\bdelta$ to be imposed to the global plate model on the inner interface $\gamma_C$ in order to generate a reaction $\lambda^C$ in balance with the reaction of the inner zone of interest submitted to the recovery of the plate displacement on its boundary $\Gamma_I$.

Starting from $\boldsymbol{\delta}_0=0$, this can be achieved through the following iterations:
\begin{enumerate}
\item\label{item:global} Run a global plate analysis with extra load $\boldsymbol{\delta}$:
\begin{equation*}
 \mathbf{u}^G_n = {\mathbf{K}^G}^{-1}\left(\mathbf{f}_{ext}^G + \boldsymbol{\delta}\right)
\end{equation*}
\item Post-process the reaction on $\gamma_C$:
\begin{equation*}
\boldsymbol{\lambda}^{C}_n = \left(\mathbf{K}^{C}\mathbf{u}^G_{n|\omega_C} - \mathbf{f}^{C}_{ext}\right)_{|\gamma_C}
\end{equation*}
\item Recover the 3D displacement on $\Gamma_I$:
\begin{equation*}
 \mathbf{u}^{L}_{n|\Gamma^I} = \mathbf{R}\mathbf{u}^{G}_{n|\Gamma_I} +\mathbf{w}^G(\mathbf{u}^{G})_n
\end{equation*}
\item\label{item:local} Solve the local problem with imposed displacement on $\Gamma_I$:
\begin{equation*}
\left\{\begin{aligned}&\mathbf{f}_{int}^L(\mathbf{u}^L_n) +\mathbf{f}_{ext}^L =0  \\
& \mathbf{u}^{L}_{n|\Gamma_I}  \text{ given} 
\end{aligned}\right.
\end{equation*}
\item Post-process the local reaction of $\Gamma_C$:
\begin{equation*}
\boldsymbol{\lambda}^{L}_n = -\left(\mathbf{f}_{int}^L(\mathbf{u}^L_n)  +\mathbf{f}_{ext}^L\right)^{\iI}_{|\Gamma_C}
\end{equation*}
\item Compute the residual on $\gamma_C$:
\begin{equation*}
\mathbf{r}_n=\boldsymbol{\lambda}^{C}_n + \mathbf{R}^{T}\boldsymbol{\lambda}^{L}_n 
\end{equation*}
\item If residual is small enough then exit, else update the extra load:
\begin{equation*}
\bdelta_{n+1|\gamma_C} = \bdelta_{n|\gamma_C} - \mathbf{r}_n
\end{equation*}
and go back to \ref{item:global}.
\end{enumerate}

One important point is that the \textit{global step} \ref{item:global} and the \textit{local step} \ref{item:local} can be processed with different software. For the test-case developed here, Code\_Aster is used for the global plate problem, and COFAST3D is used for the local 3D contact problem.

In the end the local/global method is a fixed point algorithm similar to a modified Newton algorithm. The rate of convergence of this algorithm can be quite slow. However acceleration techniques can be applied:
\begin{itemize}
\item quasi-Newton acceleration like SR1 algorithm \cite{gendre2011two} or BFGS,
\item dynamic relaxation like Aitken's algorithm ($\Delta^2$) \cite{liu2014},
\item mixed boundary conditions \cite{gendre2011two},
\item Krylov solvers (only in the linear case) \cite{guguin14phd}.
\end{itemize}

Based on \cite{guguin14phd}, it appears that, among those methods,  SR1 quasi-Newton acceleration \eqref{sec:SR1} technique leads to better performance while being simple to implement in a non-intrusive manner.

Note that running one global analysis (step 1) followed by a local reanalysis with given Dirichlet conditions (steps 3-4) without iterations corresponds to the industrialists' practice called submodeling (or sometimes structural zoom). Such an approach is ``purely descending'' in the sense that there is no feedback from the local computation towards the global scale.

\section{Analysis of the iterative corrections at the global and local levels }\label{sec:glob-local-res}

In order to evaluate the convergence of the algorithm, the relative norm of the residual at iteration $i$ is defined as:
\begin{equation}
r^{rel}_i = \frac{\Vert \boldsymbol{\lambda}^{C}_i + \mathbf{R}^{T}\boldsymbol{\lambda}^{L}_i \Vert_{2}}{\Vert \boldsymbol{\lambda}^{C}_i \Vert_{2}}
\end{equation}

In this section, the corrections associated with the coupling along the iterations are analyzed for the fixed-point algorithm (i.e. without acceleration techniques) with convergence threshold set to $10^{-6}$ which is clearly enough for mechanical quantities to have converged.

\begin{figure}[ht]
\centering
%\subfloat[][Evolution of the residual during iterations.]{\label{fig:residual}
%\includegraphics[width=0.475\linewidth]{pg_0001.pdf}
%}
%\subfloat[][Traction $F_{x}$ along $\Gamma_I^{sup}$ (top plate.)]{\label{fig:GL_conv}
%\includegraphics[width=0.45\linewidth]{pg_0002.pdf}}
\includegraphics[width=0.9\linewidth]{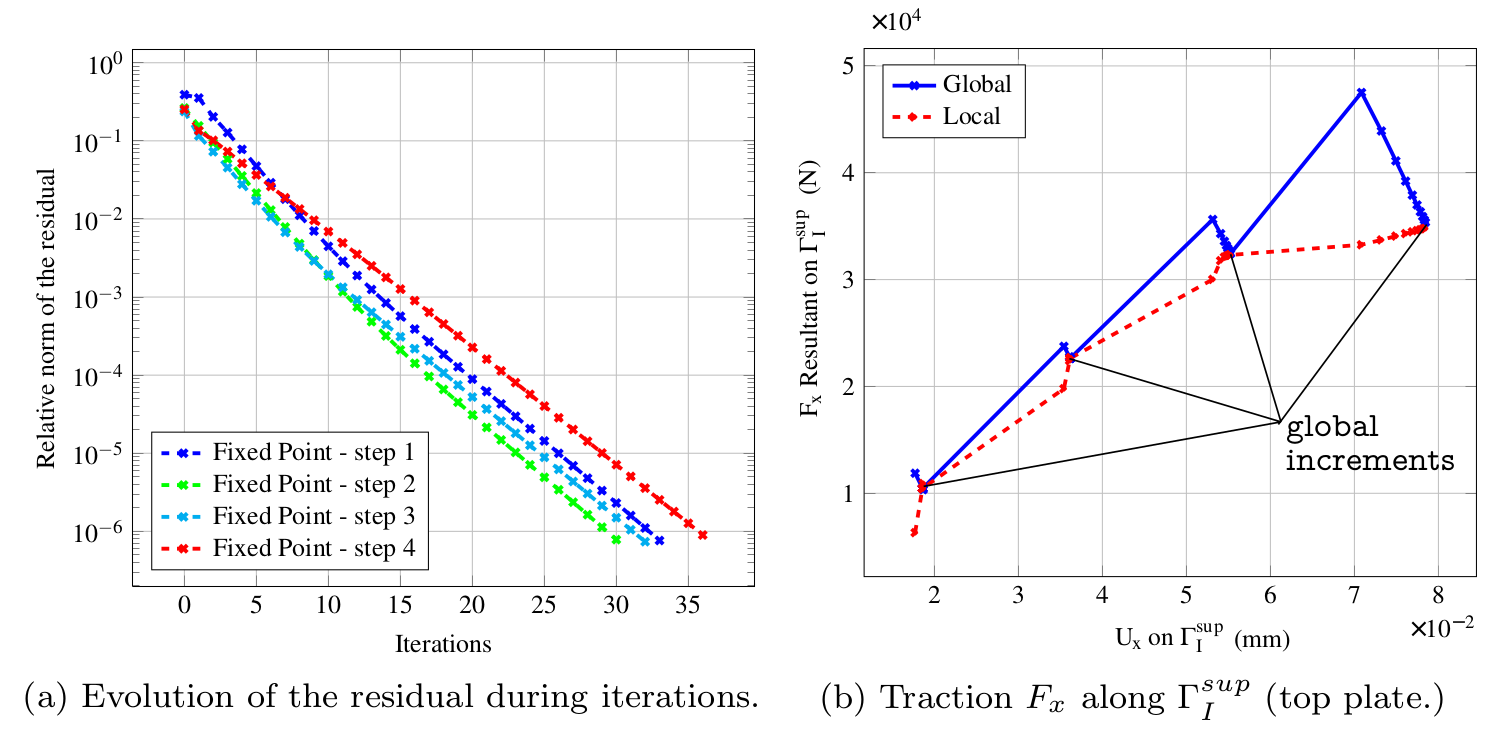}
\caption{Convergence of the fixed-point algorithm: (a) Relative norm of the residual (b) Evolution of one component of the global and local nodal reactions.}\label{fig:07}
\end{figure}

%Without using an acceleration technique 
A typical convergence curve is shown \Figure~\ref{fig:07}(a).%, note that in general SR1 acceleration allows to divide the number of iterations by a factor 3. 
This curve shows that the convergence is more or less independent of the level of nonlinearity in the bolt (weak for the first global time step, strong for the last one). This type of result is typical of nonlinear relocalization methods \cite{Cresta07,Hinojosa2014,Bordeu2009,allix2010}. Let us note that, as shown in \Figure~\ref{fig:07}(b), four global time steps are used for the plate computation, which appears to be sufficient, see Section~\ref{sec:params}. Anyhow, sub-stepping is used for the computation of the bolt; balance between the local and global model is only reached at the global time steps  in  \Figure~\ref{fig:07}(b). For more involved problems, a better control of the  global time steps and of the local sub-stepping process could be obtained using error indicators on the time discretization like in  \cite{allixkerfriden2010}.

\begin{remark}
With the proposed technique, for each global increment the first iteration corresponds to a classical submodeling approach: this enables us to easily measure the quality of a that approach, which is a question often raised by engineers. In this application, the level of the error of the submodeling approach is about $25\%$.
\end{remark}

\subsection{Global effect of the correction}\label{sec:glo-correc}

\subsubsection{Effect of the preload of the bolt}

\begin{figure}[ht]
\centering
\includegraphics[width=0.55\linewidth]{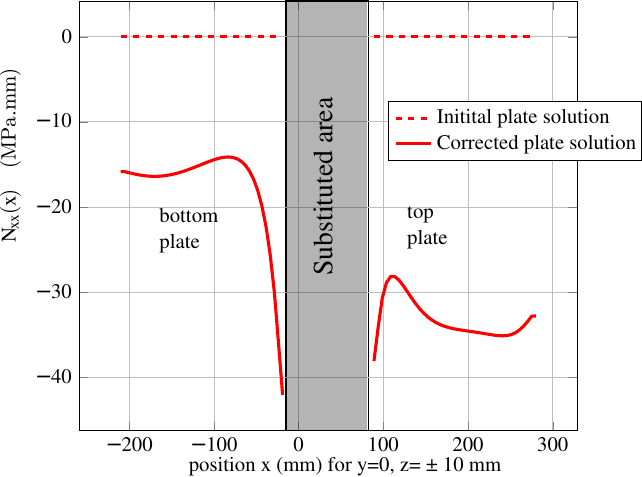}
	\caption{ Effect of the tightening (preload): value of the $\mathbf{N}_{xx}$ component of the plate stress tensor.}\label{fig:08}
\end{figure}

The global effect of the tightening of the bolt cannot be predicted using a plate theory because the associated forces are equal to zero. Nevertheless this tightening has a non-negligible global effect as can be seen in figure \Figure~\ref{fig:08}.

\subsubsection{Additional global effects of the bolt on the plate solution}

\begin{figure}[ht]
\centering
\includegraphics[width=0.7\linewidth]{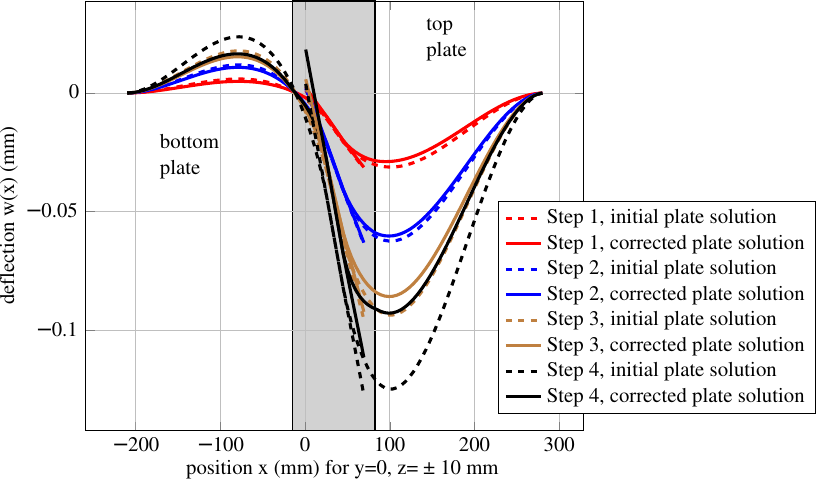}
	\caption{Deformation of the plates, before and after iterations.}\label{fig:09}
\end{figure}

After the preload of the bolt, the plate is loaded in tension, in four global time steps. The analysis of the number of global increments is presented in Section~\ref{sec:params}. 
The global solution is modified along  the iterations to match with the 3D model of the bolt, as can be seen on \Figure~\ref{fig:09} which shows the values of the transverse displacement in the mid-section of the plate for the initial plate solution and the corrected one. 

One can notice that the correction increases  with the time steps, and it becomes significant for the third step and very large for the last one. 
This is of course due to the fact that for the first two time steps the bolt acts more or less as a 3D rigid connector. Whereas the sliding of the bolt becomes significant for the third time step. During the last stage of the loading  loosening of the bolt happens, as can be checked from the local analysis of the solution presented in the next sub-section. Let us recall that the initial plate solution also corresponds to what would be obtained by a submodeling approach since submodeling only improves the local reanalyzed area.

\subsection{Local effect of the correction}\label{sec:loc_correc}

In this section the solution within the bolt is compared for three different approaches: a reference full 3D simulation, a submodeling approach and the mixed 2D-3D model with $L=0$ (minimal size of the local 3D model) obtained at the convergence of the iterations (convergence threshold is $10^{-6}$). % connection between the two models just at the limit of the bolt as computed by the non intrusive approach for a residual tolerance of $10^{-6}$ and the one corresponding to a sub-modelling approach (first iteration of the non-intrusive scheme when using a fixed point algorithm. 

\begin{figure}[ht]\centering
\includegraphics[width=0.8\linewidth]{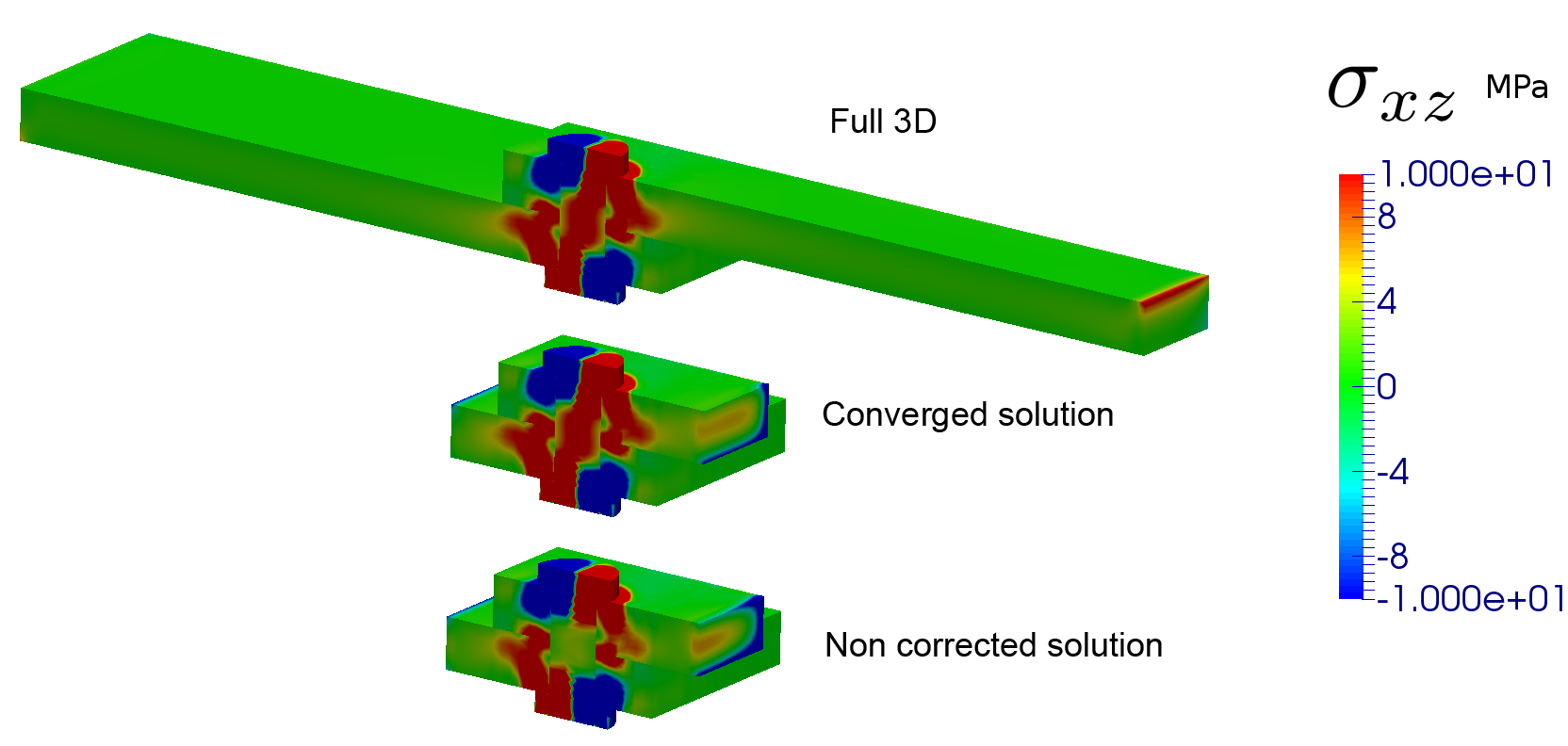}
\caption{Transverse shear stress disbribution:  reference (top), converged 2D-3D solution (middle), submodeling (bottom).}\label{fig:10}
\end{figure}

For example the figure \Figure~\ref{fig:10} shows the comparison of the $\sigma_{xz}$ stress field between the submodeling solution, the hybrid solution and the reference one. It can be seen that the solution associated with the non-intrusive approach matches very well the 3D reference solution while the submodeling gives inaccurate results  inside the screw where the shear stress is greatly underestimated.   

\begin{figure}[ht]
\centering
%\subfloat[][Response of the bolt on $\Gamma_C^{sup}$.]{\label{fig:forces_vs_depl_diff}
%\includegraphics[width=0.45\linewidth]{pg_0005.pdf}
%
%}
%\subfloat[][$N_{xx}$ vs sliding near the hole.]{\label{fig:sliding_vs_jump_diff}
%\includegraphics[width=0.45\linewidth]{pg_0006.pdf}
%
%}\\
%\subfloat[][Dissipation between the plate.]{\label{fig:dissipation}
%\includegraphics[width=0.5\linewidth]{pg_0007.pdf}
%}
%\subfloat[][Extraction points.]{\label{fig:extraction_point}
%\includegraphics[width=0.4\linewidth]{extraction3D.pdf}
%}
\includegraphics[width=0.9\linewidth]{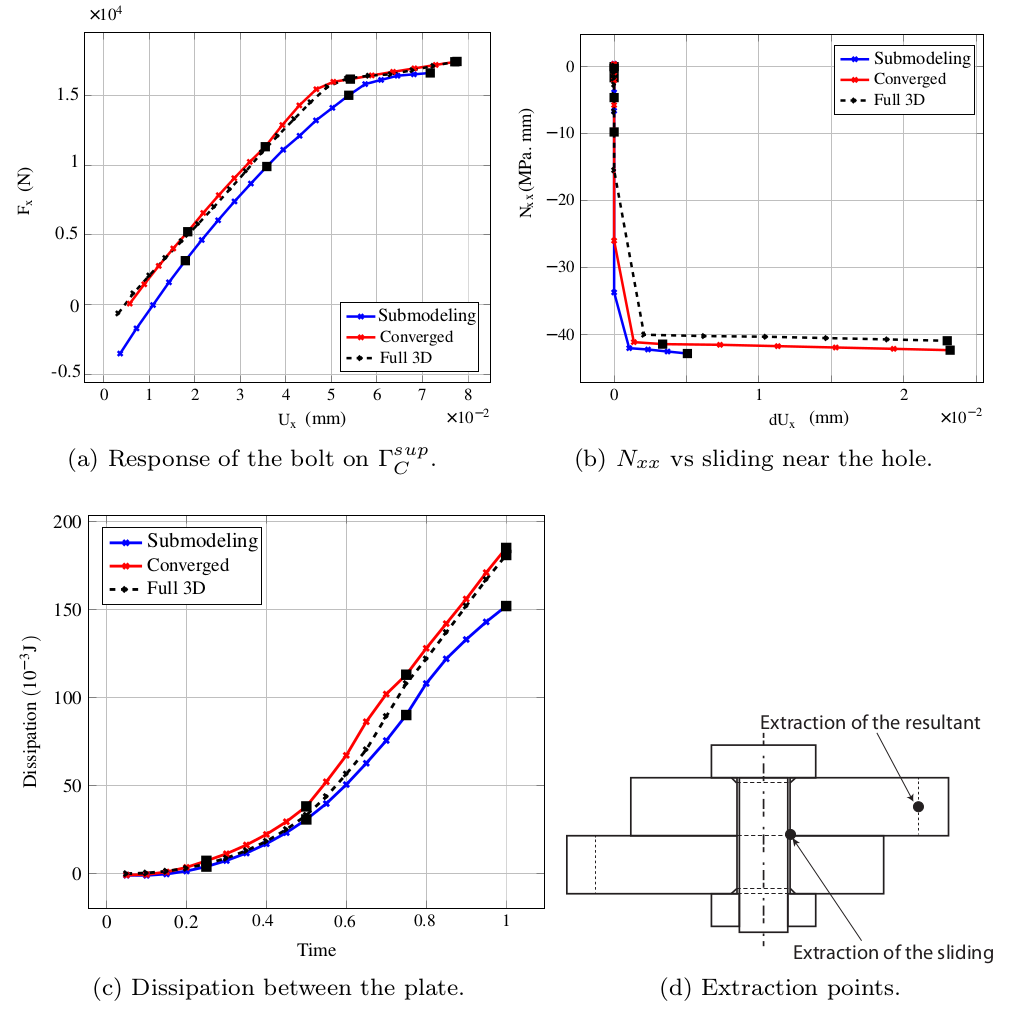}
\caption{Comparison between full 3D, submodeling and 2D-3D approaches.}\label{fig:11}
\end{figure}

We now compare global quantities associated to the bolt (like the amount of dissipation due to sliding or the global axial component of the force acting on the joint) as well as local quantities (like the pointwise values of the sliding) for the three models. On  the figures  \Figure~\ref{fig:11}(a) and \Figure~\ref{fig:11}(b), local stress quantities are extracted on a point defined on \Figure~\ref{fig:11}(d).
It appears that the relative residual provides an efficient indicator for the errors on global quantities. Typically the submodeling  leads to 20\% of error on the dissipation. Note that the error committed by submodeling on local quantities can be much larger. In comparison, the converged 2D-3D solution and the full 3D solution are quite close.

\begin{figure}[ht]
\centering
%\subfloat[][Preload.]{
%\includegraphics[width=0.4\linewidth]{glissement_0.png}
%}
%\subfloat[][Step 1.]{
%\includegraphics[width=0.4\linewidth]{glissement_1.png}
%}
%\subfloat[][Step 2.]{
%\includegraphics[width=0.4\linewidth]{glissement_2.png}
%}\\
%\subfloat[][Step 3.]{
%\includegraphics[width=0.4\linewidth]{glissement_3.png}
%}
%\subfloat[][Step 4.]{
%\includegraphics[width=0.4\linewidth]{glissement_4.png}
%}
\includegraphics[width=0.8\linewidth]{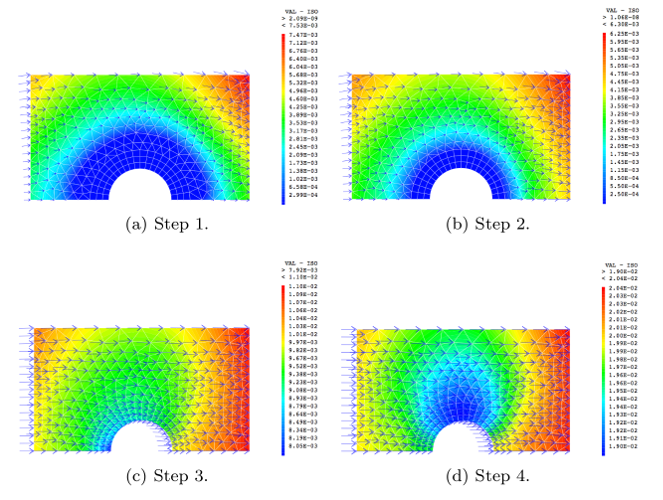}
\caption{Sliding at each global time step of the converged solution (mm).}\label{fig:12}
\end{figure}

\begin{figure}[ht]
\centering
%\subfloat[][Preload.]{
%\includegraphics[width=0.4\linewidth]{glissement_0.png}
%}
%\subfloat[][Step 1.]{
%\includegraphics[width=0.4\linewidth]{diff_GLI_1.pdf}
%}
%\subfloat[][Step 2.]{
%\includegraphics[width=0.4\linewidth]{diff_GLI_2.pdf}
%}\\
%\subfloat[][Step 3.]{
%\includegraphics[width=0.4\linewidth]{diff_GLI_3.pdf}
%}
%\subfloat[][Step 4.]{
%\includegraphics[width=0.4\linewidth]{diff_GLI_4.pdf}
%}
\includegraphics[width=0.8\linewidth]{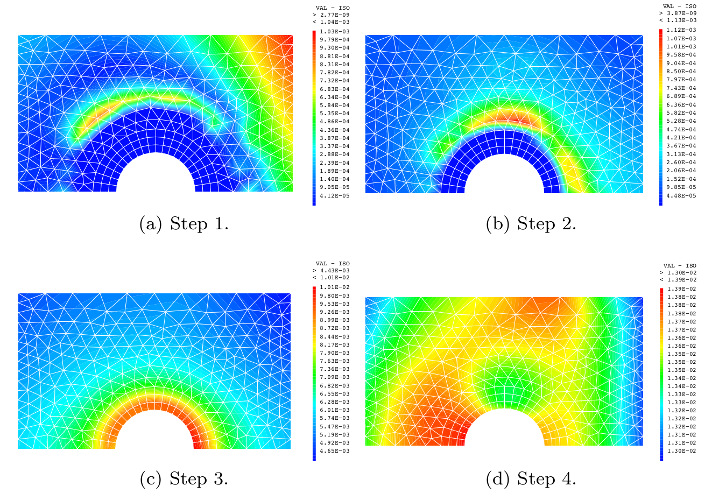}
\caption{Maps of the difference of sliding between converged and submodeling solutions at each global time step (mm).}\label{fig:13}
\end{figure}

Figure \Figure~\ref{fig:12} presents the sliding obtained by the converged 2D-3D model; it matches the 3D reference (few percents of deviation on the maximal sliding). On the contrary, the submodeling approach underestimates the sliding especially for the third and fourth global time steps, as can be seen on figure \Figure~\ref{fig:13} where the error on the maximal sliding is more than $30\%$.

To conclude, the coupling approach seems reliable contrarily to the submodeling approach for which the level of error is $20-25\%$ on global quantities and can be much more for local quantities of interest. The lack of conservatism of the submodeling extends similar results obtained in previous studies on localized plasticity or buckling for example. 

\section{Control of some parameters of the method and acceleration technique }\label{sec:params}

As already discussed several parameters can be tuned for the non-intrusive 2D-3D coupling strategy :
\begin{itemize}
\item the number of global time steps,
\item the choice of the position of the interface in the global model (dimension $L$ in \Figure~\ref{fig:06}),
\item the width of the buffer zone (dimension $b$ in \Figure~\ref{fig:06}).
\end{itemize}

The size of the buffer zone has been chosen from our experience on the 2D-3D coupling involving composite plates and orthotropic plies \cite{guguin2014nonintrusive},  in order to minimize the problem of artificial edge effects between the two models. Besides the parameters of the nonlinear solver used for the bolt computation have been tuned in order to ensure a high precision according to dedicated papers like \cite{Champaney1999249}.

\subsection{Influence of the number of global time steps}

\begin{figure}[ht]
\centering
%\subfloat[][Observation of the bolt response on $\Gamma_C^{sup}$.]{
%\includegraphics[width=0.45\linewidth]{pg_0008.pdf}
%}
%\subfloat[][$N_{xx}$ vs sliding near the screw.]{
%\includegraphics[width=0.45\linewidth]{pg_0009.pdf}
%}
\includegraphics[width=0.9\linewidth]{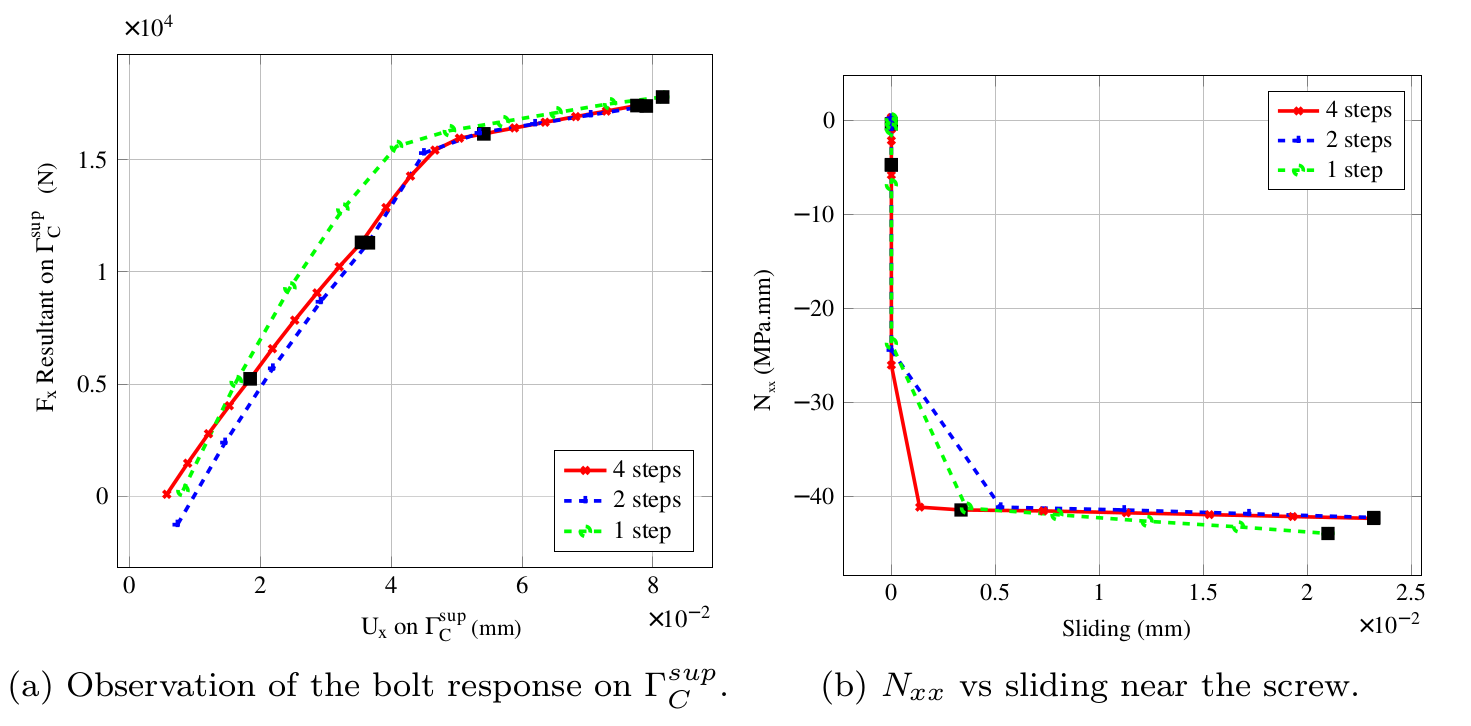}
\caption{Comparison of solutions for different number of global time steps.}\label{fig:14}
\end{figure}

The number of time steps has been simply chosen on the basis of the comparison of solutions, starting from one global time step and increasing this number. On the figure \Figure~\ref{fig:14}, a comparison between $1$, $2$, and $4$ time steps is shown.
The solutions with $2$ and $4$ time steps are considered here as quite close, whereas the one using only one time step deviates largely in the middle of the load sequence. 
Let us note that the use of large time steps at the global level is made possible thanks to the chosen algorithm with nonlinear localization. Indeed the possibility to use sub-stepping at the local scale greatly reduces the potential difficulties of convergence.

\subsection{Influence of the position of the interface between the 2D and the 3D models}

The question of the position of the interface raises in fact the question of the validity of the plate theory with respect to the 3D theory and has been largely discussed in \cite{guguin2014nonintrusive}. 
From what is known on the validity of the plate and shell theories, in the case of isotropic materials, one expects the 2D-3D model to be a good approximation of the 3D reference, for an interface situated from the bolt at a distance superior to the thickness of the plate.
It is therefore interesting to analyze the influence of  the  position of the interface and to compare the cases of $L=0$ mm and $L=15$ mm (see \Figure~\ref{fig:06}).

\begin{figure}[ht]
\centering
%\subfloat[][Residual during the iteration]{\label{fig:residual_position}
%\includegraphics[width=0.45\linewidth]{pg_0010.pdf}
%}
%\subfloat[][$N_{xx}$ vs sliding near the hole.]{\label{fig:slide_position}
%\includegraphics[width=0.45\linewidth]{pg_0011.pdf}
%
%}\\
%\subfloat[][Axial displacement.]{\label{fig:deformed_L0L15}
%\includegraphics[width=0.9\linewidth]{pg_0012.pdf}
%}
\includegraphics[width=0.9\linewidth]{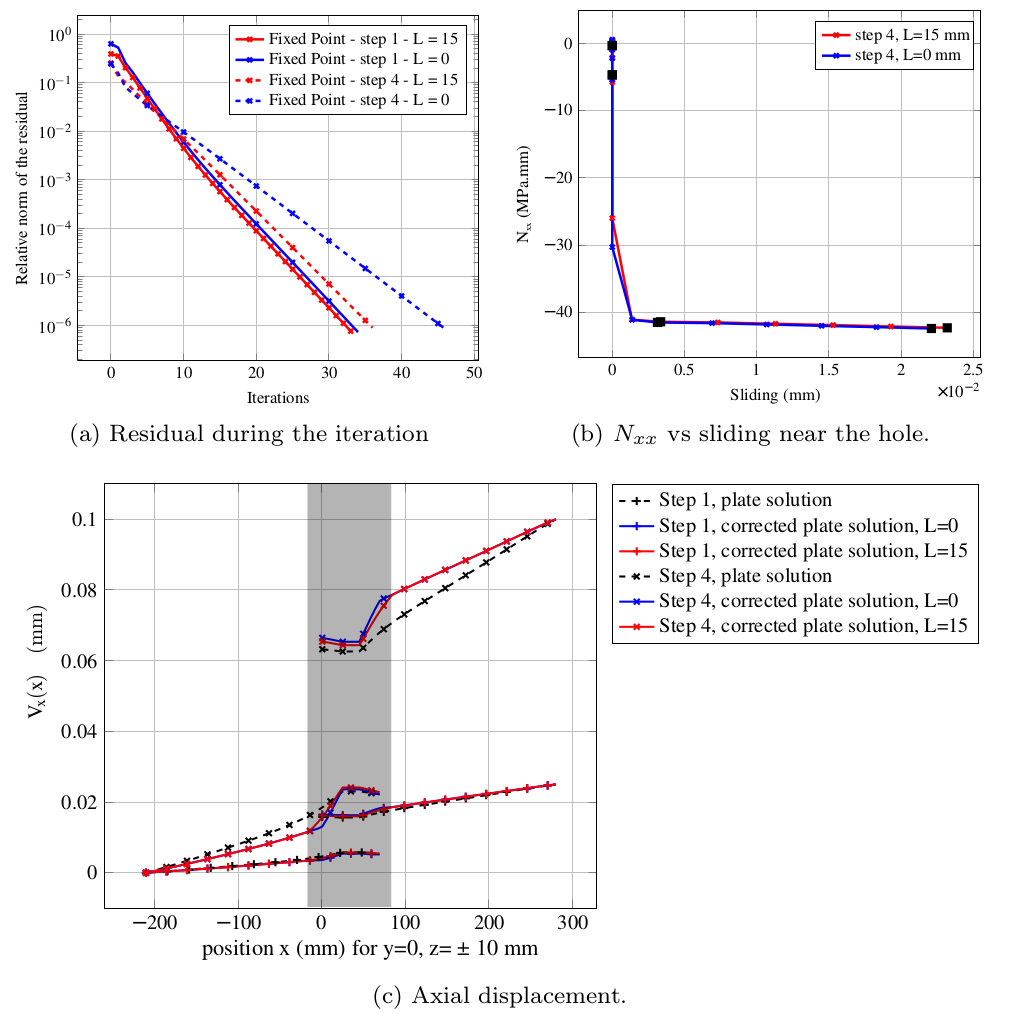}
\caption{Comparison of the results for $L=0$ mm and $L=15$ mm}\label{fig:15}
\end{figure}

From a global point of view, the final deformed shapes of the corrected solutions in \Figure~\ref{fig:15}(c) are very close in the common plate domain. They only slightly differ in the local area of interest which are different for the two models.

\begin{figure}[ht]\centering
\includegraphics[width=0.8\linewidth]{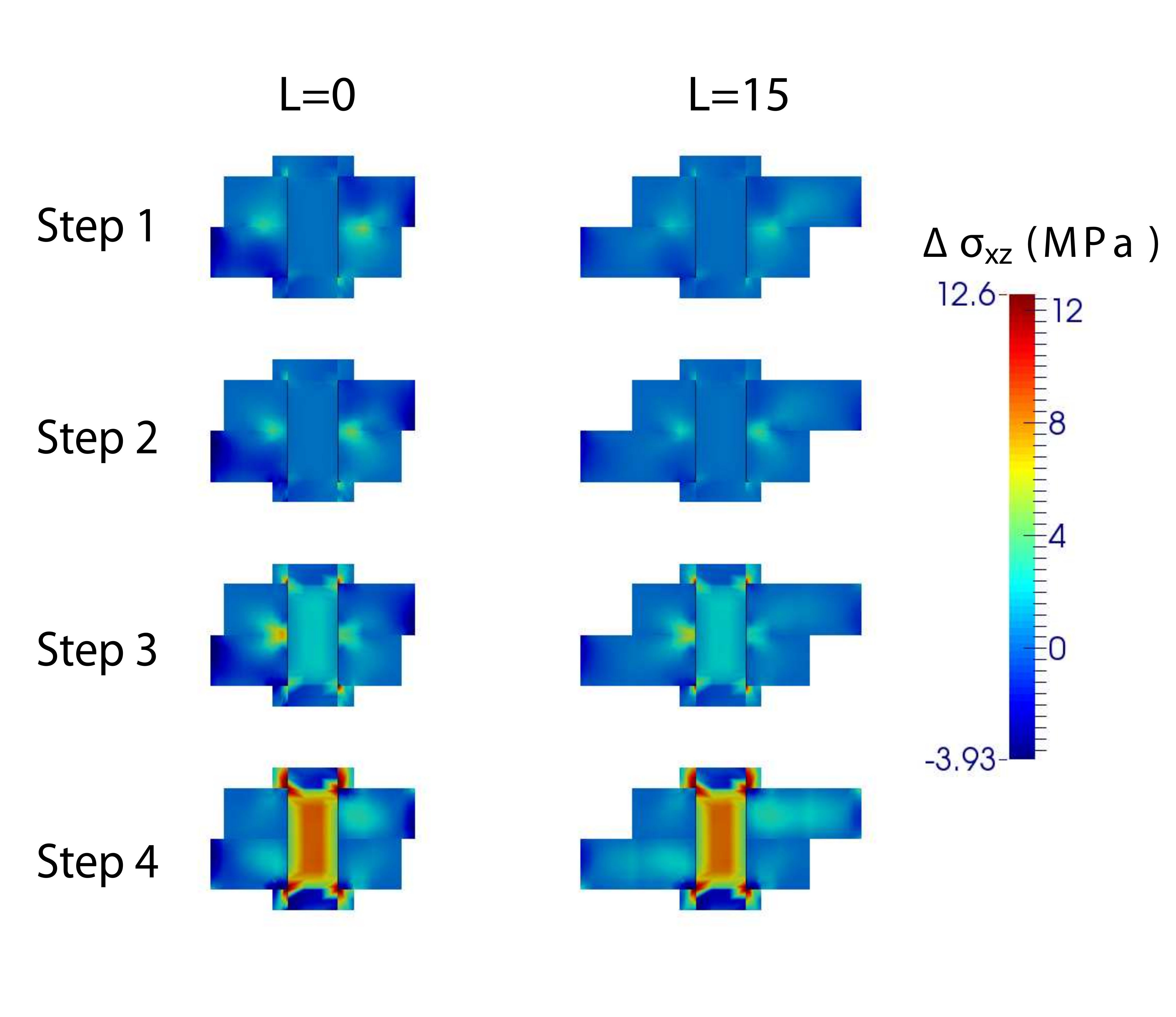}
\caption{Error in the local domain for two positions of the interface and at each time step.}\label{fig:SIXZ_erreur}
\end{figure}

\begin{figure}[ht]
\centering
\includegraphics[width=0.5\linewidth]{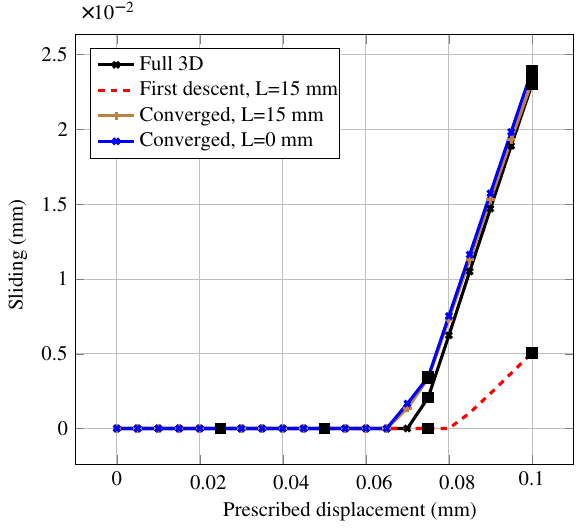}
	\caption{Sliding between the two plates near the nut for different configurations.}\label{fig:comp_position}
\end{figure}

From a  local point of view, on figure \Figure~\ref{fig:15}(b),  both cases give close results. This is confirmed when analyzing the difference of the shear stress within the bolt $\Delta \sigma_{xz} = \sigma_{xz}^{conv} - \sigma_{xz}^{0}$, as shown in the \Figure~\ref{fig:SIXZ_erreur}. 
As already analyzed, if for the first two steps, most of the differences are localized around the nut, during the sliding phase the correction mostly concerns the nut itself. 
The \Figure\ref{fig:comp_position} shows the evolution of the local tangential jump with respect to the prescribed displacement. This is an interesting quantity in order to observe the initiation of the sliding.
The 2D-3D models give close predictions for both values of $L$. They are much more closer to the reference than what is predicted by the simple submodeling.

On the figure \Figure~\ref{fig:15}(a), it appears, as expected, that a larger number of iterations has to be carried out when the interface is located on the bolt boundaries ($L=0$). This effect of the interface location on the convergence rate is corrected when using a SR1 acceleration technique as can be seen in the next subsection.

\subsection{Quasi-Newton acceleration}\label{sec:SR1}

A final component of the non-intrusive substitution is the implementation of acceleration techniques. In that context, the SR1 quasi-Newton acceleration was tested in \cite{gendre2011two}, and the Aitken Delta-2 dynamic relaxation method in \cite{duval2014non}. For this study, we only implement SR1 acceleration.

%We only tested the SR1 acceleration. This acceleration is based on the Sherman-Morrison formula, which produce a rank one modification of the operator ${\mathbb{K}^G}^{-1}$ based on the previous residual. With this correction at each step the correction operator are closer to the coupled operator sought, and the rate of convergence is largely improved.

%The rate of convergence of the algorithm then improve through a non-intrusive update of ${\mathbb{K}^G}^{-1}$ during the correction step, in order to have a correctional operator closer to ${\mathbb{K}^{hyb}}^{-1}$.
\begin{figure}[ht]
\centering
%\subfloat[][In the first step]{\label{fig:residual_buffer}
\includegraphics[width=0.5\linewidth]{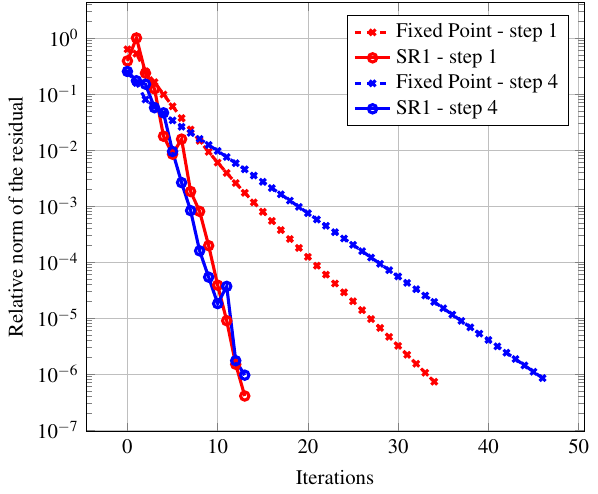}
\caption{Effect of the SR1 acceleration on the global-local rate of convergence.}\label{fig:residual_SR1}
\end{figure}

The \Figure~\ref{fig:residual_SR1} presents convergence plots. %, when the SR1 acceleration is applied, a relative $0.01\%$ error in the residual is achieved in less than 10 iterations and 
The convergence is roughly 3 times faster with SR1. Moreover the rate of convergence is almost independent from the load step and the position of the interface.
%The rate of convergence is approximately three time faster than for the fixed point algorithm. Moreover the convergence is nearly independent of the location of the interface. 
This property can be explained by the fact that the corrections induced by the SR1 acceleration technique are adapted to take into account the main differences between the 2D and 3D models. 

\section{Conclusions}

In this paper, a non-intrusive coupling between plate models and 3D models \cite{guguin2014nonintrusive} has been extended to deal with the precise computations of bolted plates: the plate model with simplified connector was coupled with a full 3D nonlinear model of the bolt. The flexibility of the method was exploited to easily define the coupled model and to use of two dedicated pieces of software:  Code\_Aster for the plate computations, and COFAST3D for the nonlinear computation of the bolt (including many surfaces of friction).

The proposed technique enabled us to analyze the possibilities and limits of the use of plate connectors and submodeling approach in that case. It appears that, even when the bolt response is globally linear, the 3D effects induced by the bolt largely modify the plate solution itself. This shows that a rigid connector is a poor model to describe the connection between two plates. As expected, when important sliding occurs, such modeling becomes irrelevant. Moreover, the submodeling technique may lead to significant local errors and non-conservative results. 

Another important feature for future applications concerning the treatment of multiple bolts in interaction, is that using an interface located at the limit of the bolt leads to acceptable results when compared to the reference solution. Works on that type of problem is in progress.

Another issue concerns the modeling of the bolt itself. In practice, some of the parameters of the bolted assembly are not precisely determined, as the preload of the nut or the friction coefficient. Such types of problems have been analyzed by dedicated techniques for the bolt computation \cite{gant2010modeling,Gant2013128}, including multiresolution \cite{vidal02}. The use of the proposed non-intrusive techniques allows to extend these studies to the case of 2D-3D structural analyses in a straightforward manner. In addition, the use of model reduction techniques for the global model itself, as proposed in \cite{kerfriden2012local} should lead to a very significant reduction of the computational time.  

%%%%%%%%%%%%%%%%
%% Background %%
%%
%\nocite{oreg,schn,pond,smith,marg,hunn,advi,koha,mouse}

%%%%%%%%%%%%%%%%%%%%%%%%%%%%%%%%%%%%%%%%%%%%%%
%%                                          %%
%% Backmatter begins here                   %%
%%                                          %%
%%%%%%%%%%%%%%%%%%%%%%%%%%%%%%%%%%%%%%%%%%%%%%

\section*{Competing interests}
  The authors declare that they have no competing interests.

%\section*{Author's contributions}
%    Text for this section \ldots

\section*{Acknowledgements}

This work was partially funded by the French National Research Agency as part of project ICARE (ANR-12-MONU-0002-04).

%%  After submission of the .TEX file,                     %%
%%  you will be prompted to submit your .BBL file.         %%

% if your bibliography is in bibtex format, use those commands:
%\bibliographystyle{bmc-mathphys} % Style BST file (bmc-mathphys, vancouver, spbasic).

\end{document}